\documentclass{article}

\usepackage{amsmath,amsthm,amssymb,bigstrut,mathrsfs,game}

\newtheorem{theorem}{Theorem}[section]
\newtheorem{definition}[theorem]{Definition}
\newtheorem{proposition}[theorem]{Proposition}

\newtheorem{lemma}[theorem]{Lemma}

\renewcommand{\L}{\ensuremath{\mathscr{L}}}
\newcommand{\N}{\ensuremath{\mathscr{N}}}
\renewcommand{\P}{\ensuremath{\mathscr{P}}}
\newcommand{\R}{\ensuremath{\mathscr{R}}}

\newcommand{\negone}{\overline{1}}
\renewcommand{\Star}{\ast}

\begin{document}

\title{Mis\`ere Canonical Forms of Partizan Games \\ {\small PREPRINT}}
\author{Aaron N. Siegel\\Institute for Advanced Study\\1 Einstein Drive\\Princeton, NJ 08540}

\maketitle

\begin{abstract}
We show that partizan games admit canonical forms in mis\`ere play.  The proof is a synthesis of the canonical form theorems for normal-play partizan games and mis\`ere-play impartial games.  It is fully constructive, and algorithms readily emerge for comparing mis\`ere games and calculating their canonical forms.

We use these techniques to show that there are precisely 256 games born by day~2, and to obtain a bound on the number of games born by day~3.
\end{abstract}

\section{Introduction}

The identification of a canonical form theorem is often a pivotal moment in understanding a particular theory of combinatorial games.  Canonical forms provide evidence of cohesive structure and reassurance that we are not floundering about in uncharted wilderness.

Among theories of finite, loopfree games in disjunctive compounds, there have been three major results in this direction.  The oldest is the celebrated \emph{Sprague--Grundy Theorem}~\cite{grundy_1939,sprague_1935,sprague_1937}: every normal-play impartial game is equivalent to a Nim-heap.  Equally important is Conway's generalization to normal-play partizan games~\cite{berlekamp_1982,conway_1976}: every such game $G$ is equivalent to a unique game $G'$ with no dominated or reversible followers.  Finally, there is an analogous theorem, also due to Conway, for mis\`ere-play impartial games~\cite{conway_1976,siegel_200Xg}.

The existence of these three theorems leaves an obvious and glaring exception.  Despite some recent progress due to Mesdal and Ottaway~\cite{mesdal_2007}, the disjunctive theory of mis\`ere-play partizan games has largely remained a mystery.  In this paper, we show that such games do indeed admit canonical forms, and that they are not all that different from normal-play canonical forms.  Furthermore, the proof is constructive, and in fact canonical forms are no harder to compute in mis\`ere play than in normal play.

\section{Mis\`ere Equivalence}

In normal play, there is a simple recursive test for equality: $G \geq H$ iff no $G^R \leq H$ and $G \leq$ no $H^L$.  Furthermore, canonical forms can be obtained by eliminating dominated options and bypassing reversible ones.

In this paper we generalize these results to mis\`ere play.  Remarkably, the definitions of dominated and reversible options, and the canonical form theorem itself, carry over to mis\`ere play without modification.  However, the recursive test for $\geq$, and the associated proofs, are considerably more involved.

We denote by $o^-(G)$ the mis\`ere outcome class of $G$:
\[
o^-(G) = \begin{cases}
\L & \textrm{if Left can win no matter who plays first}; \\
\R & \textrm{if Right can win no matter who plays first}; \\
\P & \textrm{if second player (the Previous player) can win}; \\
\N & \textrm{if first player (the Next player) can win}.
\end{cases}
\]
The four outcome classes are naturally partially ordered by ``favorability to Left,'' with $\L \geq \P \geq \R$ and $\L \geq \N \geq \R$, and this induces the usual partial order on all mis\`ere games.
\[G \geq H \textrm{ iff } o^-(G+X) \geq o^-(H+X) \textrm{ for all games } X.\]
We also define mis\`ere equality in the usual manner.
\[G = H \textrm{ iff } G \geq H \textrm{ and } H \geq G.\]

We begin with a simple, but useful, proposition.

\begin{proposition}
\label{prop:geqsimpler}
$G \geq H$ iff the following two conditions hold:
\begin{enumerate}
\item[(i)] For all $X$ with $o^-(H+X) \geq \P$, we have $o^-(G+X) \geq \P$; and
\item[(ii)] For all $X$ with $o^-(H+X) \geq \N$, we have $o^-(G+X) \geq \N$.
\end{enumerate}
\end{proposition}

\begin{proof}
$\Rightarrow$ is immediate.  For the converse, we must show that $o^-(G+X) \geq o^-(H+X)$, for all $X$.  If $o^-(H+X) = \R$, then there is nothing to prove; if $o^-(H+X) = \P$ or $\N$, it is immediate from (i) or (ii), respectively.  Finally, if $o^-(H+X) = \L$, then by (i) and (ii) we have $o^-(G+X) \geq$ both $\P$ and $\N$, whence $o^-(G+X) = \L$.
\end{proof}

\section{Ends and Adjoints}

In normal play, every ``one-sided'' game $\{G^L|~\}$ or $\{~|G^R\}$ is equivalent to an integer.  In mis\`ere play, by contrast, such games represent significant pathologies and are the source of much complication.  We will sometimes use a dot to indicate ``no moves,'' when it is useful to clarify the notation; for example, we may write $\Game{0}{\cdot}$ in place of $\Game{0}{}$.

\begin{definition}
$G$ is a \emph{Left (Right) end} if $G$ has no Left (Right) option.
\end{definition}

In normal play, it is always the case that $G + (-G)$ is equal to zero.  This is emphatically false in mis\`ere play.  Nonetheless, we can give an explicit example of a game $G^\circ$ such that $G + G^\circ$ is always a \P-position.  Readers familiar with the impartial theory will recognize it as the partizan analogue of Conway's \emph{mate}.

\begin{definition}
\label{def:adjoint}
The \emph{adjoint} of $G$, denoted $G^\circ$, is given by
\[
G^\circ =
\begin{cases}
\Star & \textrm{if  $G = 0$} \\
\Game{(G^R)^\circ}{0} & \textrm{if $G \neq 0$ and $G$ is a Left end} \\
\Game{0}{(G^L)^\circ} & \textrm{if $G \neq 0$ and $G$ is a Right end} \\
\Game{(G^R)^\circ}{(G^L)^\circ} & \textrm{otherwise}
\end{cases}
\]
\end{definition}

\begin{proposition}
\label{prop:gplusgadj}
$G + G^\circ$ is a \P-position.
\end{proposition}

\begin{proof}
By symmetry, it suffices to show that Left can win $G + G^\circ$ moving second.  By definition, $G^\circ$ is not a Right end, so Right must have a move.  If Right moves to $G^R + G^\circ$ or $G + (G^L)^\circ$, Left makes the mirror image move on the other component, which wins by induction.  Finally, if~$G$ is a Left end and Right moves to $G + 0$, then Left has no move, and so wins \emph{a priori}.
\end{proof}

\begin{theorem}
\label{theorem:dualdistinguishability}
If $G \not\geq H$, then:
\begin{enumerate}
\item[(a)] There is some $T$ such that $o^-(G+T) \leq \P$ but $o^-(H+T) \geq \N$; and
\item[(b)] There is some $U$ such that $o^-(G+U) \leq \N$ but $o^-(H+U) \geq \P$.
\end{enumerate}
\end{theorem}

\begin{proof}
We know by Proposition~\ref{prop:geqsimpler} that one of (a) or (b) must hold, so it suffices to show that (a) $\Rightarrow$ (b) and (b) $\Rightarrow$ (a).  The arguments are identical, so we will show that (a) $\Rightarrow$ (b).

Fix $T$ so that $o^-(G+T) \leq \P$ and $o^-(H+T) \geq \N$, and put
\[U = \Game{\left(H^R\right)^\circ}{T}.\]
Now from $G+U$, Right has a winning move, to $G+T$.  Therefore $o^-(G+U) \leq \N$.  Likewise, consider $H+U$.  It is certainly not a Right end, since Right has a move from $U$ to $T$.  Now if Right moves to $H^R+U$, Left has a winning response to $H^R+(H^R)^\circ$.  If instead Right moves to $H+T$, then since $o^-(H+T) \geq \N$, Left wins \emph{a priori}.  Therefore $o^-(H+U) \geq \P$, as needed.
\end{proof}

Recently Mesdal and Ottaway~\cite{mesdal_2007} showed that $G \neq 0$ unless $G$ is identically zero.  The following Lemma generalizes that result, and it proves to be a crucial piece of the analysis.

\begin{lemma}
\label{lemma:leftendnotleq}
If $H$ is a Left end and $G$ is not, then $G \not\geq H$.
\end{lemma}

\begin{proof}
Put
\[T = \Big\{\left(H^R\right)^\circ~\Big|\Big|~\cdot~\Big|~\left(G^L\right)^\circ\Big\}.\]
Consider $H+T$.  If Right moves to $H^R + T$, then Left can respond to $H^R + (H^R)^\circ$, winning; if Right moves to $H + \Game{\cdot}{(G^L)^\circ}$, then Left wins outright, since he has no further move.  Therefore $o^-(H+T) \geq \P$.

Now consider $G+T$.  Right has a move to $G + \Game{\cdot}{(G^L)^\circ}$.  Left's only response is to $G^L + \Game{\cdot}{(G^L)^\circ}$, which must exist since $G$ is not a Left end.  But Right may then respond to $G^L+(G^L)^\circ$, which wins.  Therefore $o^-(G+T) \leq \N$.

This shows that $o^-(G+T) \not\geq o^-(H+T)$, so in fact $G \not\geq H$.
\end{proof}

\section{Dominated and Reversible Options}

It is a remarkable fact that the definitions of dominated and reversible options are \emph{exactly} the same as in normal play.

\begin{definition}
Let $G$ be a game.
\begin{enumerate}
\item[(a)] A Left option $G^L$ is said to be \emph{dominated} if $G^{L'} \geq G^L$ for some other Left option $G^{L'}$.
\item[(b)] A Right option $G^R$ is said to be \emph{dominated} if $G^{R'} \leq G^R$ for some other Right option $G^{R'}$.
\item[(c)] A Left option $G^L$ is said to be \emph{reversible} if $G^{LR} \leq G$ for some Right option $G^{LR}$.
\item[(d)] A Right option $G^R$ is said to be \emph{reversible} if $G^{RL} \geq G$ for some Left option $G^{RL}$.
\end{enumerate}
\end{definition}

\begin{lemma}
Suppose $G^{L_1}$ is dominated by $G^{L_2}$, and let $G'$ be the game obtained by eliminating $G^{L_1}$ from $G$.  Then $G = G'$.
\end{lemma}

\begin{proof}
Since the Left options of $G'$ are a subset of those of $G$, and since $G'$ still has at least one Left option (namely, $G^{L_2}$), we trivially have $G' \leq G$.  Thus it suffices to show that $G' \geq G$.

So fix $X$, and suppose that Left can win $G+X$ playing first (or second).  He follows exactly the same strategy on $G'+X$, except when it recommends a move from some $G+Y$ to $G^{L_1}+Y$.  In that case, we have $o^-(G^{L_1}+Y) \geq \P$.  Since $G^{L_2} \geq G^{L_1}$, this necessarily implies $o^-(G^{L_2}+Y) \geq \P$.  So Left can win by moving from $G'+Y$ to $G^{L_2}+Y$.
\end{proof}

\begin{lemma}
Suppose $G^{L_1}$ is reversible through $G^{L_1R_1}$, and let $G'$ be the game obtained by bypassing $G^{L_1}$:
\[G' = \Game{G^{L_1R_1L},G^{L'}}{G^R},\]
where $G^{L'}$ is understood to range over all Left options of $G$ except $G^{L_1}$.  Then $G = G'$.
\end{lemma}

\begin{proof}
First suppose Left can win playing first (or second) on $G + X$.  Fix a winning strategy for Left, and assume that it recommends a move on the $X$ component unless the only winning move is on $G$.  Left follows exactly the same strategy on $G' + X$ except when it recommends a move from $G$ to $G^{L_1}$.  In that case the position must be $G' + Y$, with $o^-(G^{L_1}+Y) \geq \P$.  Thus $o^-(G^{L_1R_1}+Y) \geq \N$.  Now $G$ is not a Left end (since it has $G^{L_1}$ as an option), so by Lemma~\ref{lemma:leftendnotleq}, neither is $G^{L_1R_1}$.  Therefore Left must have a winning move from $G^{L_1R_1}+Y$.  It cannot be to $G^{L_1R_1}+Y^L$, since this would imply $o^-(G+Y^L) \geq \P$, contradicting Left's choice of strategy.  Therefore Left's move to $G^{L_1R_1L} + Y$ must be winning, and he can make this move directly from $G' + Y$.

Now suppose Right can win playing first (or second) on $G + X$.  Fix a winning strategy for Right.  She follows exactly the same strategy on $G' + X$ except when Left moves from $G'$ to $G^{L_1R_1L}$.  In that case, the position must be $G^{L_1R_1L} + Y$, with $o^-(G+Y) \leq \P$.  Since $G^{L_1R_1} \leq G$, we have $o^-(G^{L_1R_1}+Y) \leq \P$, so $G^{L_1R_1L} + Y$ is necessarily a win for Right moving first.  Finally, since $G^{L_1R_1}$ is not a Left end (as noted above), neither is $G'$.  Thus if Right follows this strategy, Left can never run out of moves prematurely.
%
%
\end{proof}

\section{The Canonical Form Theorem}

\begin{theorem}[Canonical Form Theorem]
\label{theorem:canonicalform}
Suppose $G = H$, and assume that neither $G$ nor $H$ has any dominated or reversible options.  Then for every $H^L$ there is a $G^L$ such that $G^L = H^L$, and vice versa; and likewise for Right options.
\end{theorem}

In order to prove Theorem~\ref{theorem:canonicalform}, we must generalize some machinery from Conway's proof that \emph{impartial} games admit mis\`ere canonical forms.

\begin{definition}
\label{def:linked}
~
\begin{enumerate}
\item[(a)] $G$ is \emph{downlinked to $H$ (by $T$)} iff $o^-(G+T) \leq \P$ and $o^-(H+T) \geq \P$.
\item[(b)] $G$ is \emph{uplinked to $H$ (by $T$)} iff $o^-(G+T) \geq \P$ and $o^-(H+T) \leq \P$.
\end{enumerate}
\end{definition}

\begin{theorem}
\label{theorem:geqiffnolink}
$G \geq H$ iff the following four conditions hold.
\begin{enumerate}
\item[(i)] $G$ is downlinked to no $H^L$;
\item[(ii)] No $G^R$ is downlinked to $H$;
\item[(iii)] If $H$ is a Left end, then so is $G$;
\item[(iv)] If $G$ is a Right end, then so is $H$.
\end{enumerate}
\end{theorem}

\begin{proof}
For $\Rightarrow$~(i), fix any game $T$.  If $o^-(G+T) \leq \P$, then necessarily $o^-(H+T) \leq \P$ as well.  Therefore $o^-(H^L+T) \leq \N$, so $T$ cannot downlink $G$ to $H^L$.  $\Rightarrow$~(ii) is similar, and $\Rightarrow$~(iii) and~(iv) are just restatements of Lemma~\ref{lemma:leftendnotleq} (and its mirror image).

We now prove $\Leftarrow$.  First fix $T$ such that $o^-(H+T) \geq \P$, and suppose (for contradiction) that $o^-(G+T) \leq \N$.  Then either $o^-(G+T^R) \leq \P$, or $o^-(G^R+T) \leq \P$, or else $G+T$ is a Right end.  If $o^-(G+T^R) \leq \P$, then by induction on the birthday of $T$ we may assume that $o^-(H+T^R) \leq \P$, contradicting the assumption that $o^-(H+T) \geq \P$.  If $o^-(G^R+T) \leq \P$, then $T$ downlinks $G^R$ to $H$, contradicting (ii).  Finally, if $G+T$ is a Right end, then in particular $G$ is a Right end, so by (iv) $H$ is a Right end.  Therefore $H+T$ is a Right end, contradicting the assumption that $o^-(H+T) \geq \P$.

Finally, we must show that if $o^-(H+T) \geq \N$, then $o^-(G+T) \geq \N$.  The proof is identical, with (i) and (iii) in place of (ii) and (iv).
\end{proof}

\begin{theorem}
\label{theorem:linkiffnogeq}
$G$ is downlinked to $H$ iff no $G^L \geq H$ and $G \geq$ no $H^R$.
\end{theorem}

\begin{proof}
Suppose $T$ downlinks $G$ to $H$, so that $o^-(G+T) \leq \P$ and $o^-(H+T) \geq \P$.  Then necessarily $o^-(G^L+T) \leq \N$ and $o^-(H^R+T) \geq \N$, so $T$ witnesses both $G^L \not\geq H$ and $G \not\geq H^R$.

Conversely, suppose that no $G^L \geq H$ and $G \geq$ no $H^R$.  Then for each $G^L_i$, Theorem~\ref{theorem:dualdistinguishability} yields an $X_i$ such that
\[o^-(G^L_i+X_i) \leq \P \textrm{ and } o^-(H+X_i) \geq \N.\]
Likewise, for each $H^R_j$, there is some $Y_j$ such that
\[o^-(G+Y_j) \leq \N \textrm{ and } o^-(H^R_j+Y_j) \geq \P.\]
Put
\[T = \begin{cases}
\Star & \textrm{if $G = H = 0$} \\
\Game{0}{(H^L)^\circ} & \textrm{if $G = 0$ and $H$ is a nonzero Right end} \\
\Game{(G^R)^\circ}{0} & \textrm{if $H = 0$ and $G$ is a nonzero Left end} \\
\Game{Y_j,(G^R)^\circ}{X_i,(H^L)^\circ} & \textrm{otherwise}
\end{cases}
\]
We claim that $G$ is downlinked to $H$ by $T$.  We will show that $o^-(G+T) \leq \P$; the proof that $o^-(H+T) \geq \P$ is identical.

We first show that $G+T$ has a Left option.  If $G$ has a Left option, this is automatic.  If $G$ or $H$ has a Right option, then $T$ necessarily has a Left option.  This exhausts every case except when $G = 0$ and $H$ is a Right end; but then Left's move to $0$ is built into the definition of $T$.

Thus $G+T$ is not a Left end, and it therefore suffices to show that every Left option is losing.  If Left moves to $G^L_i+T$, Right can respond to $G^L_i+X_i$, which wins by choice of $X_i$.  If Left moves to $G+(G^R)^\circ$, Right can respond to $G^R+(G^R)^\circ$, which wins by Proposition~\ref{prop:gplusgadj}.  Left's move to $G+Y_j$ loses automatically, by choice of $Y_j$.  The only remaining possibility is Left's additional move to $0$ in the first two cases of the definition of $T$.  But that move is only available when $G = 0$, so it ends the game immediately.
\end{proof}

\begin{proof}[Proof of Theorem~\ref{theorem:canonicalform}]
Fix $H^L$.  Since $G \geq H$, Theorem~\ref{theorem:geqiffnolink} implies that $G$ is not downlinked to $H^L$.  By Theorem~\ref{theorem:linkiffnogeq}, either $G^L \geq H^L$ or $G \geq H^{LR}$.  The latter would imply that $H \geq H^{LR}$, contradicting the assumption that $H$ has no reversible options.  So we must have $G^L \geq H^L$.

An identical argument, using the fact that $H \geq G$, shows that $H^{L'} \geq G^L$ for some $H^{L'}$.  Therefore \[H^{L'} \geq G^L \geq H^L.\]
Since $H$ has no dominated options, we must have $H^{L'} = H^L$, so that
\[H^{L'} = G^L = H^L.\]
The same argument suffices for the remaining cases.
\end{proof}

\section{Games Born by Day 2}

There are four games born by day $1$; and they are familiar from the normal-play theory:
\[0 = \Game{\cdot}{\cdot} \qquad \Star = \Game{0}{0} \qquad 1 = \Game{0}{\cdot} \qquad \negone = \Game{\cdot}{0}\]
Remarkably, they are pairwise incomparable.

\begin{proposition}
\label{prop:day1incomparable}
The four games $0$, $\Star$, $1$, and $\negone$ are pairwise incomparable.
\end{proposition}

\begin{proof}
Theorem~\ref{theorem:geqiffnolink}(iii) and (iv) immediately yield $\Star,1 \not\geq 0,\negone$ and $0,1 \not\geq \Star,\negone$.  Since $0$ is downlinked to $0$ (by $\Star$), (i) and (ii) furthermore imply $0 \not\geq 1$ and $\negone \not\geq 0$.

Now as a trivial consequence of Theorem~\ref{theorem:linkiffnogeq}, we have that $\negone$ is downlinked to $0$ and $0$ is downlinked to $1$.  It therefore follows from Theorem~\ref{theorem:geqiffnolink}(i) that $\negone \not\geq 1,\Star$ and from (ii) that $\negone,\Star \not\geq 1$.  This exhausts all possibilities.
\end{proof}

\begin{theorem}
There are 256 games born by day 2.
\end{theorem}

\begin{proof}
There are 16 subsets of $\{0,\Star,1,\negone\}$.  This gives 256 isomorphism types for games born by day~2, so it suffices to show that every (formal) game born by day~2 is canonical.

So fix such a game $G$.  By Proposition~\ref{prop:day1incomparable}, $G$ has no dominated options, so it suffices to show that $G$ has no reversible options.  Consider some $G^{LR}$.  Since~$G$ is born by day~2, $G^{LR}$ is born by day~$0$, so necessarily $G^{LR} = 0$.  But $G^{LR}$ is a Left end and $G$ is not (since it has $G^L$ as a Left option), so by Lemma~\ref{lemma:leftendnotleq} $G \not\geq G^{LR}$.
\end{proof}

Let $\mathbb{P}$ denote the set of games born by day~2.  We next describe the partial order structure of~$\mathbb{P}$.  Define
\[\begin{array}{lcl}
\mathbb{P}^+ & = & \{G \in \mathbb{P} : G \textrm{ is a nonzero Left end}\} \\
\mathbb{P}^- & = & \{G \in \mathbb{P} : G \textrm{ is a nonzero Right end}\} \\
\mathbb{P}^0 & = & \{G \in \mathbb{P} : G \textrm{ is not an end}\}
\end{array}\]
Then we can write $\mathbb{P}$ as a disjoint union
\[
\mathbb{P} = \mathbb{P}^+ \cup \mathbb{P}^- \cup \mathbb{P}^0 \cup \{0\}.
\]
Now let $\mathbb{B}_4$ denote the complete Boolean lattice of dimension~4.  Let $\mathbb{B}_4^+$ be the partial order obtained by removing the largest element from $\mathbb{B}_4$, and likewise delete the smallest element to obtain $\mathbb{B}_4^-$.  We will show that
\[
\label{eq:poisotypes}
\mathbb{P}^+ \cong \mathbb{B}_4^+; \qquad \mathbb{P}^- \cong \mathbb{B}_4^-; \qquad \mathbb{P}^0 \cong \mathbb{B}_4^+ \times \mathbb{B}_4^-. \tag{\dag}
\]
In order to characterize the structure of $\mathbb{P}$, it will then suffice to describe the interaction between components.

First note that there are certain trivial order relations, described by the following definition.

\begin{definition}
\label{def:trivialpo}
We say that \emph{$G$ trivially exceeds $H$}, and write $G \geq_T H$, iff:
\begin{enumerate}
\item[(i)] The Left options of $G$ form a superset of those of $H$;
\item[(ii)] The Right options of $G$ form a subset of those of $H$;
\item[(iii)] $G$ is a Left end iff $H$ is a Left end; and
\item[(iv)] $G$ is a Right end iff $H$ is a Right end.
\end{enumerate}
\end{definition}

As the terminology suggests, it is a trivial fact that if $G \geq_T H$, then necessarily $G \geq H$.  We will now show that on day~2, the converse holds with only a few exceptions.

\begin{theorem}
\label{theorem:poisusuallytrivial}
Fix $G$ and $H$ satisfying Definition~\ref{def:trivialpo}(iii) and~(iv), and assume that $G$ and $H$ are both born by day~2.  If $G \geq H$, then $G \geq_T H$.
\end{theorem}

\begin{proof}
We show that every Left option of $H$ is a Left option of $G$; the argument that every Right option of $G$ is a Right option of $H$ is identical.

So fix an $H^L$; by Theorem~\ref{theorem:geqiffnolink}, $G$ is not downlinked to $H^L$.  By Theorem~\ref{theorem:linkiffnogeq}, either $G^L \geq H^L$ for some $G^L$, or else $G \geq H^{LR}$.  Now since $H$ has the Left option $H^L$, it is not a Left end, whence by assumption neither is $G$.  Since $H$ is born by day~2, we know that every $H^{LR} = 0$, so by Lemma~\ref{lemma:leftendnotleq} we cannot have $G \geq H^{LR}$.  Therefore $G^L \geq H^L$ for some $G^L$.

But $G^L$ and $H^L$ are both born by day~1.  By Proposition~\ref{prop:day1incomparable}, this implies $G^L = H^L$.
\end{proof}

Now if $G$ and $H$ are in the same component of $\mathbb{P}$, then they necessarily satisfy Definition~\ref{def:trivialpo}(iii) and~(iv).  Therefore, \emph{on each component}, the relations $\geq$ and $\geq_T$ coincide.  But this immediately establishes (\ref{eq:poisotypes}).  For example, for the isomorphism $\mathbb{P}^+ \to \mathbb{B}_4^+$, we can regard $\mathbb{B}_4$ as the powerset lattice of $\{0,\Star,1,\negone\}$; then each $G$ maps to its set of Right options.

To complete the picture of $\mathbb{P}$, we must characterize the interaction between the four components.  We are concerned specifically with the case where $H$ is a Right end, but $G$ is not; or where $G$ is a Left end, but $H$ is not (the converses are ruled out by Lemma~\ref{lemma:leftendnotleq}).

\begin{theorem}
\label{theorem:pointerrelationships}
The ordering of $\mathbb{P}$ is generated by its restrictions to $\mathbb{P}^+$, $\mathbb{P}^-$, and $\mathbb{P}^0$, together with the following four relations and their mirror images.
\[
\{|\Star,1\} \geq 0 \qquad
\{\Star|\Star,1\} \geq \{\Star|\} \qquad
\{\negone|\Star,1\} \geq \{\negone|\} \qquad
\{\Star,\negone|\Star,1\} \geq \{\Star,\negone|\} 
\]
\end{theorem}

\begin{proof}
It is a simple matter to verify each of the four stated relations.  To prove the Theorem, we must show that no further ones are possible.

We first characterize those games that compare with $0$.  So suppose $G \geq 0$.  By Theorem~\ref{theorem:geqiffnolink} $G$ is necessarily a Left end, and furthermore no $G^R$ is downlinked to $0$.  By Theorem~\ref{theorem:linkiffnogeq}, the Right options of $G$ must therefore be a subset of $\{\Star,1\}$.  So either $G = 0$, or $G \geq_T \{|\Star,1\}$.   Games with $G \leq 0$ are characterized symmetrically.

Now suppose $G \geq H$, $H$ is a Right end, and $G$ is not.  Consider any $G^R$.  By Theorem~\ref{theorem:geqiffnolink}, $G^R$ is not downlinked to $H$.  Since $H$ is a Right end, it is necessarily the case that $G^{RL} \geq H$.  In particular, $G^R$ is not a Left end.  Furthermore, $G$ is born by day~2, so $G^{RL} = 0$ and we have $H \leq 0$.  By the previous argument, this implies that $H \leq_T \{\Star,\negone|\}$.

Therefore the Right options of $G$ form a subset of $\{\Star,1\}$, and the Left options of $H$ form a subset of $\{\Star,\negone\}$.  Furthermore, the Left options of $G$ form a superset of those of $H$, just as in the proof of Theorem~\ref{theorem:poisusuallytrivial}.  Therefore $G \geq H$ is implied by one of the three given relations, each representing one possibility for the Left options of $H$.
\end{proof}

\subsection*{Antichains by Day 2}

Since $\mathbb{P}$ has such a clean structure, we can get a tight bound on the number of antichains.  From a standard reference (such as~\cite[A000372]{sloane}), we find that $\mathbb{B}_4$ has 168 antichains.  This shows that $\mathbb{P}^+$ (and hence $\mathbb{P}^-$ as well) has precisely 167, since there is a unique antichain containing the largest element of $\mathbb{B}_4$.

Consider $\mathbb{P}^0$.  We have $\mathbb{B}_4 \times \mathbb{B}_4 \cong \mathbb{B}_8$, trivially.  Again from a standard reference, we find that $\mathbb{B}_8$ has 56130437228687557907788 antichains.  Since every antichain of $\mathbb{P}^0$ is an antichain of $\mathbb{B}_8$, this gives an upper bound for the number of antichains of~$\mathbb{P}^0$.

Finally, $\{0\}$ trivially admits just two antichains, $\emptyset$ and $\{0\}$.  Since every antichain on $\mathbb{P}$ restricts to an antichain on each component, this gives an upper bound of
\[M = 2 \times 167 \times 167 \times 56130437228687557907788\]
for the number of antichains of $\mathbb{P}$.  Thus we obtain an upper bound of $M^2$ games born by day~3.  This number is large, indeed (roughly $2^{183}$), but it is much smaller than $2^{512}$, the number of nonisomorphic game trees of height~3.


\section{Relationships to Other Theories}

In this section we consider how the partizan mis\`ere theory relates to other theories of combinatorial games.  Denote by $o^+(G)$ the normal-play outcome class of $G$, and define
\[\begin{array}{lcl}
G \geq^+ H & \textrm{iff} & o^+(G+X) \geq o^+(H+X) \textrm{ for every game $X$}; \bigstrut \\
G \geq^- H & \textrm{iff} & o^-(G+X) \geq o^-(H+X) \textrm{ for every game $X$}.
\end{array}
\]
$\geq^-$ is the relation that we have been calling $\geq$; for this section only, we include the minus sign for clarity.  $\geq^+$ is, of course, the usual Berlekamp--Conway--Guy inequality for normal-play partizan games.  The following result shows that $\geq^+$ is a coarsening of $\geq^-$.

\begin{theorem}
If $G \geq^- H$, then $G \geq^+ H$.
\end{theorem}

\begin{proof}
We must show that Left can get the last move in $G - H$.  Suppose Right plays to $G^R-H$ (the argument is the same if she plays to $G - H^L$).  Since $G \geq^- H$, Theorem~\ref{theorem:geqiffnolink} implies that $G^R$ is not downlinked to $H$.  By Theorem~\ref{theorem:linkiffnogeq}, either $G^{RL} \geq^- H$ or $G^R \geq^- H^R$.  By induction, we may assume that either $G^{RL} \geq^+ H$ or $G^R \geq^+ H^R$.  In the first case, Left wins by moving to $G^{RL}-H$; in the second, by moving to $G^R-H^R$.
\end{proof}

In addition to the usual equivalences $=^+$ and $=^-$, obtained by symmetrizing $\geq^+$ and $\geq^-$, we have two further equivalences when $G$ and $H$ are impartial.
\[\begin{array}{lcl}
G =^+_I H & \textrm{iff} & o^+(G+X) = o^+(H+X) \textrm{ for every impartial game $X$}; \bigstrut \\
G =^-_I H & \textrm{iff} & o^-(G+X) = o^-(H+X) \textrm{ for every impartial game $X$}.
\end{array}
\]
It is a well-known fact that $=^+_I$ is just the restriction of $=^+$ to impartial games.  It is worth pointing out that the analogous statement does \emph{not} hold for mis\`ere games.

\begin{proposition}
There exist impartial games $G$ and $H$ such that $G =^-_I H$ but $G \neq^- H$.
\end{proposition}

\begin{proof}
It is well-known that $\Star + \Star =^-_I 0$ (see \cite{conway_1976}).  However, $\Star + \Star \neq^- 0$, by Lemma~\ref{lemma:leftendnotleq}.
\end{proof}

Therefore $=^-_I$ is a strict coarsening of $=^-$.  This highlights an interesting difference between normal and mis\`ere play: there exist impartial games $G$ and~$H$ that are distinct in partizan mis\`ere play, but that are not distinguishable by any impartial game.  This behavior arises in other theories, as well; for example, there exist impartial loopy games $G$ and~$H$ that are distinct (in normal play), but that are not distinguishable by any impartial loopy game~\cite{siegel_2006a}.  Indeed, the coincidence of $=^+$ and $=^+_I$ appears to be an artifact of the special nature of short games in normal play: it is the exception rather than the rule.

\section{Partizan Mis\`ere Quotients}

Recently Thane Plambeck~\cite{plambeck_2005} observed that, if $\mathscr{A}$ is any set of impartial games, then its mis\`ere-play structure can often be simplified by localizing the mis\`ere equivalence relation to~$\mathscr{A}$.  Plambeck showed that many important aspects of the theory can be generalized to the local setting, and the structure theory of such quotients has been explored in detail; see~\cite{siegel_200Xd,siegel_200Xe}.

It is not our intention to replicate that analysis here, but merely to remark that a partizan generalization exists.  The construction is exactly the same, but instead of a bipartite monoid $(\mathcal{Q},\mathcal{P})$, we now have a \emph{tetrapartite monoid} $(\mathcal{Q},\Pi)$, where $\Pi : \mathcal{Q} \to \{\L,\R,\P,\N\}$ is the \emph{outcome partition} for $\mathcal{Q}$.

Intriguingly, such monoids have an induced partial order structure, given by
\[x \geq y \textrm{~~iff~~} \Pi(xz) \geq \Pi(yz) \textrm{ for all } z \in \mathcal{Q}.\]
If $G \geq H$, then it is certainly true that $\Phi(G) \geq \Phi(H)$.  However, the quotient may also gain new order-relations that are not present in the universe of games.  We include one example to illustrate the rich possibilities.  In the previous section we remarked that $1$ and $\negone$ are incomparable, and we have also seen that $1 + \negone \neq 0$.  In $\mathcal{Q}(1,\negone)$, however, the expected inequalities hold:
\[\Phi(\negone) > \Phi(0) > \Phi(1); \quad \textrm{and} \quad \Phi(\negone)\Phi(1) = \Phi(0);\]
and indeed we have $\mathcal{Q}(1,\negone) \cong \mathbb{Z}$, equipped with the usual partial-order structure.  We leave it to the reader to verify these assertions.

\bibliography{games}
\end{document}